\documentclass[11pt]{article}
\usepackage{epic,latexsym,amssymb,color,graphicx}
\usepackage{tikz}

\newtheorem{thm}{Theorem}
\newtheorem{lem}[thm]{Lemma}

\newtheorem{fact}{Fact}

\newtheorem{prop}{Proposition}
\newtheorem{conj}{Conjecture}

\newcommand{\qed}{$\Box$}

\newcommand{\dg}{{\rm deg}}
\newcommand{\barq}{\overline{q}}
\newcommand{\barv}{\overline{v}}

\newcommand{\barl}{\overline{\ell}}

\newcommand{\cC}{\mathcal{C}}

\newcommand{\cH}{{\cal H}}

\newtheorem{definition}[fact]{Definition}

\newcommand{\proof}{\noindent\textbf{Proof. }}
\newcommand{\2}{ \vspace{0.2cm} }
\newcommand{\1}{ \vspace{0.1cm} }

\newcommand{\ClaimX}[1]{\noindent\textbf{Claim #1:}}
\newcommand{\ClaimPF}[1]{\noindent\textbf{Proof of Claim #1.}}
\newcommand{\ClaimQED}{ \hfill {\tiny \mbox{($\Box$)}}}

\newcommand{\hyperedgetwo}[6]{
	\pgfmathsetmacro\Done{sqrt((#4-#1)^2+(#5-#2)^2)}
	\pgfmathsetmacro\angleone{(#2>#5)*(180+asin((#4-#1)/ \Done)-asin((#1-#4)/ \Done))+asin((#1-#4)/ \Done)+asin((#3-#6)/\Done)}

\pgfmathsetmacro\angleone{\angleone-360*(\angleone>0)-360*(\angleone>360)}
	\draw ([shift=(\angleone:#3)] #1,#2)--([shift=(\angleone:#6)]#4,#5);
	\pgfmathsetmacro\Dtwo{sqrt((#1-#4)^2+(#2-#5)^2)}
	\pgfmathsetmacro\angletwo{(#5>#2)*(180+asin((#1-#4)/ \Dtwo)-asin((#4-#1)/ \Dtwo))+asin((#4-#1)/ \Dtwo)+asin((#6-#3)/\Dtwo)}
	\pgfmathsetmacro\angletwo{\angletwo-360*(\angletwo>0)-360*(\angletwo>360)}
	\draw ([shift=(\angletwo:#6)] #4,#5)--([shift=(\angletwo:#3)]#1,#2);
	\draw (#1,#2)+(\angletwo:#3) arc(\angletwo:\angleone+360*(\angleone<\angletwo):#3);
	\draw (#4,#5)+(\angleone:#6) arc(\angleone:\angletwo+360*(\angletwo<\angleone):#6);
}

\begin{document}

\title{Every $4$-regular $4$-uniform hypergraph has a \\ $2$-coloring with a free vertex}

\author{Michael A. Henning${}^{1,}$\thanks{Research
supported in part by the South African
National Research Foundation and the University of Johannesburg} \, and Anders Yeo${}^{1,2}$\\
\\
${}^1$Department of Pure and Applied Mathematics\\
University of Johannesburg \\
Auckland Park, 2006 South Africa \\
\small\tt Email:  mahenning@uj.ac.za \\
\\
${}^2$Department of Mathematics and Computer Science \\
University of Southern Denmark \\
Campusvej 55, 5230 Odense M, Denmark \\
\small\tt Email:  andersyeo@gmail.com
}

\date{}
\maketitle

\begin{abstract}
In this paper, we continue the study of $2$-colorings in hypergraphs. A hypergraph is $2$-colorable if there is a $2$-coloring of the vertices with no monochromatic hyperedge. It is known (see Thomassen [J. Amer. Math. Soc. 5 (1992), 217--229]) that every $4$-uniform $4$-regular hypergraph is $2$-colorable. Our main result in this paper is a strengthening of this result. For this purpose, we define a vertex in a hypergraph $H$ to be a free vertex in $H$ if we can $2$-color $V(H) \setminus \{v\}$ such that every hyperedge in $H$ contains vertices of both colors (where $v$ has no color). We prove that every $4$-uniform $4$-regular hypergraph has a free vertex. This proves a known conjecture. Our proofs use a new result on not-all-equal $3$-SAT which is also proved in this paper and is of interest in its own right.
\end{abstract}

{\small \textbf{Keywords:} Hypergraphs; Bipartite; $2$-Colorable; Transversal; Free vertex; NAE-$3$-SAT. } \\
\indent {\small \textbf{AMS subject classification:} 05C69}

\newpage
\section{Introduction}

In this paper, we continue the study of $2$-colorings in hypergraphs. We adopt the notation and terminology from~\cite{HeYe13,HeYe15}. A \emph{hypergraph} $H = (V,E)$ is a finite set
$V = V(H)$ of elements, called \emph{vertices}, together with a finite multiset $E = E(H)$ of arbitrary subsets of $V$, called \emph{hyperedges} or simply \emph{edges}. A $k$-edge in $H$ is an edge of size~$k$ in $H$. The hypergraph $H$ is $k$-\emph{uniform} if every edge of $H$ is a $k$-edge.
The \emph{degree} of a vertex $v$ in $H$, denoted $d_H(v)$ or simply by $d(v)$ if $H$ is clear from the context, is the number of edges of $H$ which contain $v$.
The hypergraph $H$ is $k$-\emph{regular} if every vertex has degree~$k$ in~$H$. For $k \ge 2$, let $\cH_k$ denote the class of all $k$-uniform
$k$-regular hypergraphs. The class $\cH_k$ has been widely studied, both in the context of solving problems on total domination as well as in its own right, see
for example \cite{AlBru88, HeYe13, HeYe15, HeYe16, ThYe07}.

A hypergraph $H$ is $2$-\emph{colorable} if there is a $2$-coloring of the vertices with no monochromatic hyperedge.  Equivalently, $H$ is $2$-colorable if it is \emph{bipartite}; that is, its vertex set can be partitioned into two sets such that every hyperedge intersects both partite sets.  Alon and Bregman~\cite{AlBru88} established the following result.

\begin{thm} \label{Alon}
{\rm (Alon, Bregman~\cite{AlBru88})} Every hypergraph in $\cH_k$ is
$2$-colorable, provided $k \ge 8$.
\end{thm}

Thomassen~\cite{Th92} showed that the Alon-Bregman result in Theorem~\ref{Alon} holds for all $k \ge 4$.

\begin{thm} \label{Extension}
{\rm (Thomassen~\cite{Th92})} Every hypergraph in $\cH_k$ is $2$-colorable, provided $k \ge 4$.
\end{thm}

As remarked by Alon and Bregman~\cite{AlBru88} the result is not true when $k = 3$, as may be seen by considering the Fano plane.  Sufficient conditions for the existence of a $2$-coloring in $k$-uniform hypergraphs are given, for example, by Radhakrishnan and Srinivasan~\cite{RaSr00} and Vishwanathan~\cite{Vi03}. For related results, see the papers by Alon and Tarsi~\cite{AlTa92}, Seymour~\cite{Se74} and Thomassen~\cite{Th83}.

A set $X$ of vertices in a hypergraph $H$ is a \emph{free set} in $H$ if we can $2$-color $V(H)
\setminus X$ such that every edge in $H$ contains vertices of both colors (where the vertices in $X$ are not colored). A vertex is a \emph{free vertex} in $H$ if we can $2$-color $V(H) \setminus \{v\}$ such that every hyperedge in $H$ contains vertices of both colors (where $v$ has no color). In~\cite{HeYe15} it is conjectured that every hypergraph $H \in \cH_k$, with $k \geq 4$, has a free set of size $k-3$. Further, if the conjecture is true, then the bound $k-3$ cannot be improved for any $k \ge 4$, due to the complete $k$-uniform hypergraph of order $k+1$, as such a hypergraph needs two vertices of each color to ensure every edge has vertices of both colors. The conjecture is proved to hold for $k \in \{5,6,7,8\}$.  The case when $k=4$ turned out to be more difficult than the cases when $k \in \{5,6,7,8\}$ and was conjectured separately in~\cite{HeYe15}.

\begin{conj}{\rm (\cite{HeYe15})}
\label{conj1}
Every $4$-regular $4$-uniform hypergraph contains a free vertex.
\end{conj}

\section{Main Result}

Our immediate aim is to prove Conjecture~\ref{conj1}. That is, we prove the following result, which is a strengthening of the result of Theorem~\ref{Extension} in the case when $k = 4$.

\begin{thm}  \label{t:main1}
Every $4$-regular $4$-uniform hypergraph contains a free vertex.
\end{thm}

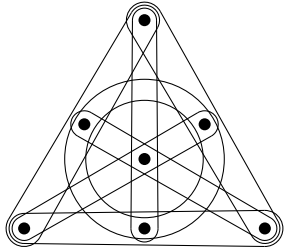
\begin{figure}[htb]
\begin{center}
\begin{tikzpicture}[scale=0.40]

\begin{scope}[xshift=0cm,yshift=0cm]
\fill (0,0) circle (0.2cm); \fill (4,0) circle (0.2cm); \fill (8,0)
circle (0.2cm); \fill (2,3.464) circle (0.2cm); \fill (6,3.464)
circle (0.2cm); \fill (4,6.928) circle (0.2cm); \fill (4,2.3093)
circle (0.2cm); \hyperedgetwo{0}{0}{0.5}{8}{0}{0.6};
\hyperedgetwo{4}{6.928}{0.5}{0}{0}{0.6};
\hyperedgetwo{4}{6.928}{0.6}{8}{0}{0.5};
\hyperedgetwo{0}{0}{0.4}{6}{3.464}{0.45};
\hyperedgetwo{4}{6.928}{0.4}{4}{0}{0.45};
\hyperedgetwo{8}{0}{0.4}{2}{3.464}{0.45}; \draw (4,2.3093) circle
(1.9523cm); \draw (4,2.3093) circle (2.6523cm);
\end{scope}

\end{tikzpicture}
\end{center}
\vskip -0.3 cm \caption{The Fano plane $F_7$} \label{f:Fano}
\end{figure}

As remarked earlier, the complete $4$-regular $4$-uniform hypergraph on five vertices has only one free vertex, and so the result of Theorem~\ref{t:main1} cannot be improved in the sense that there exist $4$-regular $4$-uniform hypergraphs with no free set of size~$2$. Theorem~\ref{t:main1} is also best possible by considering the complement, $\overline{F_7}$, of the Fano plane $F_7$, where the Fano plane is shown in Figure~\ref{f:Fano} and where its complement $\overline{F_7}$ is the hypergraph on the same vertex set $V(F_7)$ and where $e$ is a hyperedge in the complement  if and only if $V(F_7) \setminus e$ is a hyperedge in $F_7$.

Our proof of Theorem~\ref{t:main1} presented in Section~\ref{S:proof1} uses a surprising connection with not-all-equal $3$-SAT (NAE-$3$-SAT).
We will later prove a result on when NAE-$3$-SAT is not only satisfiable, but is satisfiable without assigning all variables truth values. This result is of interest in its own right, but requires some further terminology (see Section~\ref{defns}) before describing it in detail. We remark that our resulting NAE-$3$-SAT result, given by Theorem~\ref{main_nae_3sat}, has also been used by the authors in~\cite{HeYe16+} to solve a conjecture on the so-called fractional disjoint transversal number (which we do not define here). This serves as added motivation of the importance of the NAE-$3$-SAT result which can be used to solve several seemingly unrelated hypergraph problems that seem difficult to solve using a purely hypergraph approach.

\section{Terminology and Definitions}
\label{defns}

For an edge $e$ in a hypergraph $H$, we denote by $H - e$ the hypergraph obtained from $H$ by deleting the edge $e$. Two vertices $x$ and $y$ of $H$ are \emph{adjacent} if there is an
edge $e$ of $H$ such that $\{x,y\} \subseteq e$. Further, $x$ and $y$ are \emph{connected} if there is a sequence $x=v_0,v_1,v_2\ldots,v_k=y$ of vertices of $H$ in which $v_{i-1}$ is
adjacent to $v_i$ for $i=1,2,\ldots,k$. A \emph{connected hypergraph} is a hypergraph in which every pair of vertices are connected. A \emph{component} of a hypergraph $H$ is a maximal connected subhypergraph of $H$. In particular, we note that a component of $H$ is by definition connected.

A subset $T$ of vertices in a hypergraph $H$ is a \emph{transversal} in $H$ if $T$ has a nonempty intersection with every edge of $H$. In the language of transversals, a vertex $v$ is a free vertex in a hypergraph $H$ if $H$ contains two vertex disjoint transversals, neither of which contain the vertex $v$. Transversals in $4$-uniform hypergraphs are well studied (see, for example,~\cite{HeYe15,LaCh90,ThYe07}).

In order to prove Conjecture~\ref{conj1}, we use a surprising connection between an instance of not-all-equal $3$-SAT (NAE-$3$-SAT) and a $3$-uniform hypergraph. In order to state this connection we require some further terminology.

\begin{definition}
An instance, $I$, of $3$-SAT contains a set of variables, $V(I)$, and a set of clauses, $C(I)$. Each clause contains exactly three literals, which are either a variable, $v \in V(I)$, or the negation of a variable, $\barv$, where $v \in V(I)$. A clause, $c \in C(I)$, is satisfied if one of the literals in it is true. That is, the clause $c$ is satisfied if $v \in V(I)$ belongs to $c$ and $v=True$ or $\barv$ belongs to $c$ and $v=False$.  The instance $I$ is satisfied if there is a truth assignment to the variables such that all clauses are satisfied.
\end{definition}

\begin{definition}
An instance of NAE-$3$-SAT is equivalent to $3$-SAT, except that we require all clauses to contain a false literal as well as a true one.
A clause that contains both a true and false literal we call \emph{nae}-\emph{satisfiable}. If every clause in the instance $I$ is nae-satisfiable, we say that $I$ is \emph{nae}-\emph{satisfiable}.
\end{definition}

We furthermore need the following definitions.

\begin{definition}
Given an instance $I$ of NAE-$3$-SAT, we define the \textbf{associated graph} $G_I$ to be the graph with vertex set $V(I)$ and where an edge joins two variables in $G_I$ if they (either in
negated or unnegated form) appear in the same clause in $I$.

Let $I$ be an instance of NAE-$3$-SAT.  We call the instance $I$ \textbf{connected} if one cannot partition the variables $V(I)$ into non-empty sets $V_1$ and $V_2$ such that no clause contains variables from $V_1$ and $V_2$.  In other words, the graph $G_I$ associated with $I$  is connected.

A \textbf{component} of a NAE-$3$-SAT instance $I$ is a maximal connected sub-instance of $I$.
That is, the components of $I$ correspond precisely to the components of the graph $G_I$ associated with $I$.

A variable, $v \in V(I)$, is \textbf{free} if $I$ is nae-satisfiable even if we do not assign any truth value to $v$. That is, every clause in $I$ contains a true and a false literal, even without considering literals involving $v$.

The \textbf{degree} of a variable $v \in V(I)$, is the number of clauses containing $v$ or $\barv$, and is denoted by $\dg_I(v)$. If the instance $I$ is clear from the context, we simply write $\dg(v)$ rather than $\dg_I(v)$.
\end{definition}

We are now in a position to define a connection between an instance of NAE-$3$-SAT and a $3$-uniform hypergraph as follows.

\begin{definition} \label{defn2}
  If $H$ is a $3$-uniform hypergraph, we create a NAE-$3$-SAT instance $I_H$ as follows.
  Let $V(I_H) = V(H)$ and for each edge $e \in H$ add a clause to $I_H$ with the same vertices/variables in non-negated form. We call $I_H$ the NAE-$3$-SAT instance corresponding to $H$.
  Note that the instance $I_H$ is nae-satisfiable if and only if $H$ is bipartite. In fact the partite sets in the bipartition correspond to the truth values true and false.
\end{definition}

Throughout this paper, we use the standard notation $[k] = \{1,2,\ldots,k\}$.

\section{NAE-$3$-SAT}
\label{S:NAE}

In this section, we present a key result that we need in order to prove Conjecture~\ref{conj1},
namely the following theorem that establishes a fundamental property of NAE-$3$-SAT in the case when the number of clauses is less than the number of variables. An instance of NAE-$3$-SAT is non-trivial if it contains at least one variable.

\begin{thm} \label{main_nae_3sat}
Let $I$ be a connected non-trivial instance of NAE-$3$-SAT. If $|C(I)|<|V(I)|$ and $\dg_I(v) \le 3$ for all $v \in V(I)$, then $I$ is nae-satisfiable and contains a free variable.
\end{thm}
\proof Suppose, to the contrary, that the theorem is false and let $I$ be a counterexample of the theorem with minimum possible $|C(I)|$. Let $C = C(I)$ and $V = V(I)$. If $|C|=0$, then $|V|=1$ and the theorem holds, and so  $|C| \ge 1$. We will now show a number of claims which we will use to obtain a contradiction to $I$ being a counterexample.

\1

\ClaimX{A} $\dg(v) \ge 2$ for all $v \in V$.

\ClaimPF{A} Let $v \in V$ be arbitrary. If $\dg(v)=0$, then $I$ is not connected as $|C| \ge 1$, a contradiction.

Suppose that $\dg(v)=1$, and let $c \in C$ be the clause containing $v$. If $c$ contains three variables of degree~$1$, then $|C|=1$ and $|V|=3$ as $I$ is connected, and $I$ is clearly nae-satisfiable and all variables in $V$ are free, contradicting the fact that $I$ is a counterexample. Hence, the clause $c$ has at most two vertices of degree~$1$. Let $c$ contain the variables $v$, $x_1$ and $x_2$ and let $I'$ be the instance of NAE-$3$-SAT obtained by deleting $v$ and the clause $c$.

If $I'$ is connected, then, by the minimality of $|C|$, the instance $I'$ is nae-satisfiable and contains a free variable.  Assigning $v$ a truth value such that the literal containing $v$ in the clause $c$ is of opposite value to the literal containing $x_1$ or $x_2$, the instance $I$ is nae-satisfiable and contains a free variable, namely the same free variable that belongs to the instance $I'$. This contradicts the fact that $I$ is a counterexample. Therefore, $I'$ is not connected.

Since $I$ is connected and $I'$ is not connected, the instance $I'$ contains two components, one containing the variable $x_1$ and the other the variable $x_2$ that belonged to the clause $c$. Both $x_1$ and $x_2$ have degree at most~$2$ in $I'$, while the degree of all other variables in $I'$ is at most~$3$. Therefore, both components of $I'$ satisfy the condition of the theorem and, by the minimality of $|C|$, are nae-satisfiable. (We remark that there exists a free variable in both components, but in this case we assign every variable a truth value.) If the literals associated with $x_1$ and $x_2$ in the clause $c$ have the same truth values, then we can reverse the truth value of all variables in one of the components. Hence, we may assume that the literals associated with $x_1$ and $x_2$ in the clause $c$ have different truth values, implying that $I$ is nae-satisfiable and contains $v$ as a free variable. Once again, we contradict the fact that $I$ is a counterexample.

As $v$ was chosen arbitrarily we have proven Claim~A.~\ClaimQED{}

\2

\ClaimX{B} There exists a $v \in V$, such that $\dg(v) = 2$.

\ClaimPF{B} If the claim was false, then by Claim~A we would have $\dg(v)=3$ for all $v \in V$, which would imply that $3|C| = \sum_{v \in V} \dg(v) = 3|V|$, which is
a contradiction to $|C|<|V|$.~\ClaimQED{}

\2

By Claim~B, there exists a variable $v \in V$ such that $\dg(v)=2$. Let $c_1$ and $c_2$ be the clauses containing the variable $v$ and let $Q$ contain all
variables belonging to $c_1$ or $c_2$.

\2

\ClaimX{C} $|Q| \le 4$.

\ClaimPF{C} Suppose, to the contrary, that $|Q| \ge 5$. Since the clauses $c_1$ and $c_2$ both contain three variables, and $v$ belongs to both clauses, we note that $|Q| \le 5$. Consequently, $|Q| = 5$.
Let $I'$ be the NAE-$3$-SAT obtained from $I$ by deleting $c_1$, $c_2$ and $v$.

Suppose that $I'$ contains four distinct components. In this case, each component of $I'$ contains a variable from $Q \setminus \{v\}$. Possibly, a component of $I'$ may contain only one variable and no clause. By the minimality of $|C|$, the instance $I'$ is nae-satisfiable. (We remark that there exists a free variable in each of the four components, but in this case we assign every variable a truth value.) We can set the variables in $Q \setminus \{v\}$ such that both $c_1$ and $c_2$ are
nae-satisfiable by reversing all truth values in any of the components of $I'$, if required, implying that $I$ is nae-satisfiable and contains $v$ as a free variable.  This contradicts the fact that $I$ is a counterexample. Therefore, $I'$ contains at most three distinct components.

Let $Q \setminus \{v\} = \{q_1,q_2,q_3,q_4\}$. Renaming variables, if necessary, we may assume that
$q_4$ and one of $q_1$, $q_2$ or $q_3$ belong to the same component in $I'$. Renaming $q_1$, $q_2$ and $q_3$, if necessary, we may assume that $q_1$ and $q_2$ are variables in $c_1$ and $q_3$ is a variable in $c_2$. If $v$ is negated in $c_i$, then negate all literals in $c_i$, for $i \in [2]$. This does not change the problem and implies that we may assume, without loss of generality, that $v$ is not negated in both $c_1$ and $c_2$.

Let $\ell_i$ be the literal containing $q_i$ in $c_1$ for $i \in [2]$, and so $\ell_i \in \{q_i,\barq_i\}$. Further, let $\ell_i$ be the literal containing $q_i$ in $c_2$ for $i \in \{3,4\}$, and so $\ell_i \in \{q_i,\barq_i\}$. Let $c'$ be a new clause $\{\ell_1,\ell_2,\barl_3\}$ and let $I'' = I' \cup c'$. We note that $I''$ is connected and the degree of all vertices in $I''$ is at most~$3$. Further, $|C(I'')| = |C(I')|+1 = |C|-1 < |V|-1 = |V(I'')|$. In particular, $|C(I'')| < |C|$. By the minimality of $|C|$, the instance $I''$ is nae-satisfiable and contains a free variable, $f$.

If $f \not\in \{q_1,q_2,q_3\}$, then we can always assign to the variable $v$ a truth value, such that $c_1$ and $c_2$ are nae-satisfiable, as indicated in Table~1, where $T$ denotes $True$ and $F$ denotes $False$ and the cases when $(\ell_1,\ell_2,\ell_3) \in \{(T,T,F),(F,F,T)\}$ are impossible due to the clause~$c'$. Further, the free variable, $f$, in $I''$ is also a free variable in $I$. Therefore, the instance $I$ is nae-satisfiable and contains a free variable. This contradicts the fact that $I$ is a counterexample. Therefore, $f \in \{q_1,q_2,q_3\}$.

\begin{center}
\begin{tabular}{|c|c|c|c|c|c|c|c|c|c|c|} \cline{1-2} \cline{4-5} \cline{7-8} \cline{10-11}
$(\ell_1,\ell_2,\ell_3)$ & $v$ & \hspace{0.5cm}  & $(\ell_1,\ell_2,\ell_3)$ & $v$ & \hspace{0.5cm}  & $(\ell_1,\ell_2,\ell_3)$ & $v$ & \hspace{0.5cm}  & $(\ell_1,\ell_2,\ell_3)$ & $v$ \\  \cline{1-2} \cline{4-5} \cline{7-8} \cline{10-11}
$(T,T,T)$       & $F$ &                 & $(T,F,T)$       & $F$ &                 & $(F,T,T)$       & $F$ &                 & $(F,F,T)$       & N/A  \\   
$(T,T,F)$       & N/A &                 & $(T,F,F)$       & $T$ &                 & $(F,T,F)$       & $T$ &                 & $(F,F,F)$       & $T$  \\   \cline{1-2} \cline{4-5} \cline{7-8} \cline{10-11}
\end{tabular}
\begin{center}
\textbf{Table~1.} Possible assignments of truth values.
\end{center}
\end{center}

Suppose that $f \in \{q_1,q_2\}$. Renaming $q_1$ and $q_2$, if necessary, we may assume that $f=q_1$. As $c'$ is nae-satisfiable in $I''$, the literals $\ell_2$ and $\ell_3$ must have the same truth value. We can therefore assign $v$ the opposite truth value to these two literals in order to get a nae-satisfiable assignment of $I$ where the variable $q_1$ is free. This contradicts the fact that $I$ is a counterexample. Therefore, $f=q_3$.

As $c'$ is nae-satisfiable in $I''$ and $f$ is free in $I''$, the literals $\ell_2$ and $\ell_3$ must have opposite truth values. We can therefore assign $v$ the opposite truth value to the literal $\ell_4$ in order to get a nae-satisfiable assignment of $I$ where the variable $q_3$ is free. Once again, this contradicts the fact that $I$ is a counterexample.
This completes all cases and therefore completes the proof of Claim~C.~\ClaimQED{}

\2

\ClaimX{D} $|Q| \le 3$.

\ClaimPF{D} Suppose, to the contrary, that $|Q| \ge 4$, which by Claim~C implies that $|Q|=4$.  Therefore, there must exist variables $q_1$, $q_2$ and $q_3$ such that $c_1$ contains the variables $v$, $q_1$ and $q_2$ and the clause $c_2$ contains $v$, $q_1$ and $q_3$. As in the proof of Claim~C.3, we may assume, without loss of generality, that $v$ is not negated in both $c_1$ and $c_2$. Let $I'$  be the NAE-$3$-SAT obtained from $I$ by deleting the two clauses $c_1$ and $c_2$, and deleting the variable $v$.

Let $\ell_i$ be the literal containing $q_i$ in $c_1$ for $i \in [2]$, and so $\ell_i \in \{q_i,\barq_i\}$. Further, let $\ell_3$ be the literal containing $q_3$ in $c_2$, and so $\ell_3 \in \{q_3,\barq_3\}$. Let $c'$ be a new clause $\{\ell_1,\ell_2,\barl_3\}$ and let $I'' = I' \cup c'$. We note that $I''$ is connected and the degree of all vertices in $I''$ is at most~$3$. Further, $|C(I'')| = |C(I')|+1 = |C|-1 < |V|-1 = |V(I'')|$. In particular, $|C(I'')| < |C|$. By the minimality of $|C|$, the instance $I''$ is nae-satisfiable and contains a free variable, $f$.

If $f \not\in \{q_1,q_2,q_3\}$, then proceeding exactly as in the proof of Claim~C, we show that the instance $I$ is nae-satisfiable and contains a free variable, a contradiction. Therefore, $f \in \{q_1,q_2,q_3\}$.

If $f=q_1$, then as $c'$ is nae-satisfiable in $I''$, the literals $\ell_2$ and $\ell_3$ must have the same truth value. We can therefore assign $v$ the opposite truth value to these two literals in order to get a nae-satisfiable assignment of $I$ where the variable $q_1$ is free.

If $f=q_2$, then as $c'$ is nae-satisfiable in $I''$, the literals $\ell_1$ and $\ell_3$ must have the same truth value. We can therefore assign $v$ the opposite truth value to these two literals in order to get a nae-satisfiable assignment of $I$ where the variable $q_2$ is free.

If $f=q_3$, then as $c'$ is nae-satisfiable in $I''$, the literals $\ell_1$ and $\ell_2$ must have opposite truth values. We can therefore assign $v$ the opposite truth value to the literal corresponding to $q_1$ in $c_2$ in order to get a nae-satisfiable assignment of $I$ where the variable $q_3$ is free.

In all the above three cases, the instance $I$ is nae-satisfiable and contains a free variable, a contradiction. This completes all cases and therefore also the proof of Claim~D.~\ClaimQED{}


\medskip
By Claim~D, $|Q| \le 3$. As every clause contains three variables, this implies that $|Q|=3$. Let $Q = \{v,q_1,q_2\}$. Let $I^*$ be the an instance of NAE-$3$-SAT with $V(I^*)=\{v,q_1,q_2\}$ and $C(I^*)=\{c_1,c_2\}$.

\2

\ClaimX{E} The instance $I^*$ is nae-satisfiable and has a free variable.

\ClaimPF{E}  If at most one literal in $c_1$ is identical to those in $c_2$, then we simply reverse all literals in $c_1$. This does not change the problem and now there are at least two literals in $c_1$ that are identical with those in $c_2$. Renaming the variables, if necessary, we may assume that the literal containing $q_i$ in $c_1$ and $c_2$ for $i \in [2]$ is identical. The variable $v$ is therefore a free variable as may be seen by assigning opposite truth value to the literal containing $q_1$ and $q_2$ in $c_1$ (and therefore also in $c_2$). Thus, $I^*$ is nae-satisfiable and has a free variable.~\ClaimQED{}

\medskip
By Claim~E, the instance $I^*$ is nae-satisfiable and has a free variable. Therefore, $I \ne I^*$, implying that at least one of $q_1$ and $q_2$ has degree~$3$ in $I$. There is therefore a clause $c_3$, different from $c_1$ and $c_2$, containing $q_1$ or $q_2$. Renaming $q_1$ and $q_2$, if necessary, we may assume that $c_3$ contains $q_2$.

\2
\ClaimX{F} The clause $c_3$ does not contain the variable~$q_1$.

\ClaimPF{F}
Suppose, to the contrary, that $c_3$ contains~$q_1$. Let $q_3$ be the variable in $c_3$ which is different from $q_1$ and $q_2$. Let $I''$ be obtained from $I$ by deleting the three clauses $c_1$, $c_2$ and $c_3$, and deleting the three variables $v$, $q_1$ and $q_2$. We note that $I''$ is connected and the degree of all vertices in $I''$ is at most~$3$. Further, $|C(I'')| = |C|-3 < |V|-3 = |V(I'')|$. In particular, $|C(I'')| < |C|$. By the minimality of $|C|$, the instance $I''$ is nae-satisfiable and contains a free variable. (We remark, however, that here we do not need the fact that there exists a free variable in this case.) By Claim~E, it is possible to assign truth values to two of the variables in $\{v,q_1,q_2\}$ with the third vertex a free variable such that $c_1$ and $c_2$ are both nae-satisfiable. At least one of $q_1$ or $q_2$ has been assigned a truth value, say $q_1$. Let $\ell_1$ be the literal containing $q_1$ in $c_3$, and let $\ell_3$ be the literal containing $q_3$ in $c_3$. By reversing all truth values in $I''$, if necessary, we can guarantee that the literals $\ell_1$ and $\ell_3$ have opposite truth values. Therefore, $I$ is nae-satisfiable and one of the vertices in $\{v,q_1,q_2\}$ is free.  This is a contradiction to $I$ being a counterexample.~\ClaimQED{}

\medskip
By Claim~F, the clause $c_3$ does not contain the variable~$q_1$. Let $I'$  be the NAE-$3$-SAT obtained from $I$ by deleting the two clauses $c_1$ and $c_2$, and the variable $v$. As in the proof of Claim~E, we may assume that there are at least two literals in $c_1$ that are identical with those in $c_2$.

\2

\ClaimX{G} The variable $v$ is not free in $I^*$.

\ClaimPF{G}
Suppose, to the contrary, that $v$ is free in $I^*$. If the literal containing $q_1$ in $c_1$ and $c_2$ is not identical, then in order for the variable $v$ to be free in $I^*$, the literal containing $q_2$ in $c_1$ and $c_2$ is not identical. This contradicts our assumption that at least two literals in $c_1$ are identical with those in $c_2$. Hence, the literal containing $q_i$ in $c_1$ and $c_2$ for $i \in [2]$ is identical.

If in $c_1$ and $c_2$ exactly one of $q_1$ and $q_2$ is negated, then let $c_3'$ be a new clause obtained from $c_3$ by replacing the literal $q_2$ with $q_1$ or $\barq_2$ with $\barq_1$. If in $c_1$ and $c_2$ either both or none of $q_1$ and $q_2$ are negated, then let $c_3'$ be a new clause obtained from $c_3$ by replacing the literal $q_2$ with $\barq_1$ or $\barq_2$ with $q_1$. Let $I''$ be the instance obtained from $I'$ by deleting the clause $c_3$ and the variable $q_2$, and adding the clause $c_3'$. We note that $I''$ is connected and the degree of all vertices in $I''$ is at most~$3$. Further, $|C(I'')| = |C|-2 < |V|-2 = |V(I'')|$. In particular, $|C(I'')| < |C|$.

By the minimality of $|C|$, the instance $I''$ is nae-satisfiable and contains a free variable. (We remark, however, that here we do not need the fact that there exists a free variable in $I''$ in this case.) If in $c_1$ and $c_2$ exactly one of $q_1$ and $q_2$ is negated, then we let $q_2$ have the same truth value as $q_1$; otherwise, we let $q_2$ have the opposite truth value of $q_1$. This implies that $I$ is nae-satisfiable even without assigning $v$ a value. Therefore, the variable $v$ is free and $I$ is not a counterexample to the theorem, a contradiction.~\ClaimQED{}

\medskip
By Claim~G, the variable $v$ is not free in $I^*$. By our earlier assumptions, there are at least two literals in $c_1$ that are identical with those in $c_2$. If the literal containing $q_i$ in $c_1$ and $c_2$ is identical for $i \in [2]$, then $v$ is free in $I^*$, a contradiction. This implies that if the literal containing $q_i$ in $c_1$ and $c_2$ is not identical, then the literal containing $q_{3-i}$ in $c_1$ and $c_2$ is identical for $i \in [2]$. Further, the literal containing $v$ in $c_1$ and $c_2$ is identical for $i \in [2]$.

Let $c_3'$ be a new clause obtained from $c_3$ by replacing the literal $q_2$ with $q_1$ or $\barq_2$ with $\barq_1$. Let $I''$ be the instance obtained from $I'$ by deleting the clause $c_3$ and the variable $q_2$, and adding the clause $c_3'$. We note that $I''$ is connected and the degree of all vertices in $I''$ is at most~$3$. Further, $|C(I'')| = |C|-2 < |V|-2 = |V(I'')|$. In particular, $|C(I'')| < |C|$. By the minimality of $|C|$, the instance $I''$ is nae-satisfiable and contains a free variable.

Suppose that $q_1$ is free in $I''$.  In this case, we can assign values to $v$ and $q_2$ if $q_1$ is free in $I^*$ and to $v$ and $q_1$ if $q_2$
is free in $I^*$, such that $c_1$ and $c_2$ are nae-satisfiable. Therefore, $I$ is nae-satisfiable and contains $q_1$ or $q_2$ as a free variable. This is a contradiction to $I$ being a counterexample. Hence, $q_1$ is not free in $I''$.

Since $I''$ contains a free variable and $q_1$ is not free in $I''$, some variable, $w$, different from $q_1$ is a free variable in $I''$. We now assign $q_2$ the same truth value as $q_1$ and we assign $v$ the opposite truth value to the literal corresponding to $q_1$ in $c_1$ (or $c_2$). With this truth assignment, both $c_1$ and $c_2$ are nae-satisfiable, noting that the literal containing $v$ in $c_1$ and $c_2$ is identical for $i \in [2]$. Therefore, $I$ is nae-satisfiable and contains $w$ as a free variable. Once again, this is a contradiction to $I$ being a counterexample, which completes the proof of Theorem~\ref{main_nae_3sat}.~\qed

\medskip
We remark that the result of Theorem~\ref{main_nae_3sat} is best possible in the following sense.

\begin{prop}
For any $s \ge 1$, there exists a non-trivial connected instance $I$ of NAE-$3$-SAT with $3s$ variables satisfying $0 \le |C(I)| < |V(I))|$ and $\dg_I(v) \le 3$ for all $v \in V(I)$ such that $I$ is nae-satisfiable and contains exactly one free variable.
\end{prop}
\proof Let $s \ge 1$ and let $I$ be an instance of NAE-$3$-SAT with variables $V(I) = \{ v_i^j \mid i \in [s] \mbox{ and } j \in [3]\}$ and clauses $C(I) = C_1 \cup C_2$, where
\[
\begin{array}{rcl}
C_1 & = & \{ (v_i^1,v_i^2,v_i^3\}, (\barv_i^1, v_i^2, v_i^3) \mid i \in [s] \} \1 \\
 C_2 & = & \{ (v_i^1, v_{i+1}^2, \barv_{i+1}^3) \mid i \in [s-1] \}. \\
\end{array}
\]

By construction, $I$ is connected and the degree of all vertices in $I$ is at most~$3$. Further, $|C(I)| = 3s-1 = |V(I)|-1$. We will now show that the only free variable in $I$ is $v_s^1$. Due to the clauses in $C_1$ we note that $v_i^2$ and $v_i^3$ must be assigned opposite truth values in any nae-satisfiable truth assignment for all $i \in [s]$. For every $i \in [s-1]$, we note that $v_{i+1}^2$ and $\barv_{i+1}^3$ have the same truth value since $v_{i+1}^2$ and $v_{i+1}^3$ have opposite truth values. This implies that $v_i^1$ must be assigned the opposite truth value to $v_{i+1}^2$ and $\barv_{i+1}^3$. This is true for all $i \in [s-1]$, which implies that $v_s^1$ is the only variable that can be free. It is not difficult to see that $v_s^1$ is free and this also follows from Theorem~\ref{main_nae_3sat}, noting that none of the other variables are free.~\qed

\section{Proof of Theorem~\ref{t:main1}}
\label{S:proof1}

Using Theorem~\ref{main_nae_3sat}, we prove Theorem~\ref{t:main1}. First, we present the following lemma.

\begin{lem} \label{lem1}
Let $H$ be a connected $3$-uniform hypergraph with no isolated vertex. If $H$ has fewer edges than vertices and has maximum degree at most~$3$, then $H$ contains at least two free vertices.
\end{lem}
\proof
Let $H$ be a connected $3$-uniform hypergraph with no isolated vertex. Suppose that $H$ has fewer edges than vertices and has maximum degree at most~$3$. Let $I_H$ be the
NAE-$3$-SAT instance corresponding to $H$.  By Theorem~\ref{main_nae_3sat}, the instance $I_H$ is nae-satisfiable and has a free vertex, say $v$. Assigning color~$1$ to true variables and color~$2$ to false variables we obtain a $2$-coloring $\cC$ of $H$ where $v$ has no color and all hyperedges of $H$ contain vertices of both colors. Let $E_v$ be all edges in $H$ containing the vertex $v$. Since $H$ has no isolated vertex, we note that $|E_v| = d_H(v) \ge 1$.

We say that a vertex in $H$ is $\cC$-fixed if in some edge in $E(H) \setminus E_v$ it is the only vertex of its color in the $2$-coloring $\cC$. We note that every edge of $H$ is a $3$-edge, and every edge in $E(H) \setminus E_v$ contains vertices of both colors in $\cC$. Thus, in every edge in $E(H) \setminus E_v$ there is a vertex whose color is unique in that edge. Thus, every edge in $E(H) \setminus E_v$ gives rise to exactly one vertex that is $\cC$-fixed. Therefore, there are at most $|E(H) \setminus E_v|$ vertices in $H$ that are $\cC$-fixed.

By supposition, $|E(H)| < |V(H)|$. Hence,
$|E(H) \setminus E_v| \le |E(H)| - |E_v| \le
(|V(H)|-1) - d_H(v) \le |V(H)| - 2$, implying that at least two vertices in $H$ are not $\cC$-fixed.
Clearly, the vertex $v$ is not $\cC$-fixed. Let  $u$ be a vertex different from $v$ that is not $\cC$-fixed. Renaming colors if necessary, we may assume that $u$ has color~$1$. Thus, every edge in $E(H) \setminus E_v$ that contains $u$ contains another vertex of color~$1$ and one vertex of color~$2$. Let $\cC'$ be the coloring obtained from $\cC$ by removing the color~$1$ from $u$ and assigning color~$1$ to $v$. Since $\cC$ is a $2$-coloring of $H$, so too is $\cC'$ a $2$-coloring of $H$. However, in the $2$-coloring $\cC'$ the vertex $u$ is a free vertex. Thus, $H$ has at least two free vertices, namely $u$ and $v$.~\qed

\medskip
We are now in a position to prove Theorem~\ref{t:main1}. Recall its statement.

\medskip
\noindent \textbf{Theorem~\ref{t:main1}}. \emph{Every $4$-regular $4$-uniform hypergraph contains a free vertex.} \\

\noindent \textbf{Proof of Theorem~\ref{t:main1}.} We may assume that $H$ is connected as otherwise we consider each component of $H$ separately. By Thomassen's Theorem~\ref{Extension}, there exists a $2$-coloring, $\cC$, of $H$ such that no edge of $H$ is monochromatic. Analogously to the proof of Lemma~\ref{lem1}, we call a vertex $\cC$-fixed if in some edge in $E(H)$ it is the only vertex of its color in $\cC$.  If some vertex is not $\cC$-fixed, then it is a free vertex, as we can simply uncolor it. Therefore, we may assume that all vertices in $H$ are $\cC$-fixed, for otherwise the desire result follows.

For every edge $e \in E(H)$, let $v^*(e)$ be the vertex of unique color in $e$, if such a vertex exists. By assumption, all vertices in $H$ are $\cC$-fixed, implying that for every vertex $u \in V(H)$ we have $u = v^*(e)$ for some edge $e$ in $H$. Since $H$ is a $4$-regular $4$-uniform hypergraph, we note that $|V(H)|=|E(H)|$. Thus since every vertex in $H$ is $\cC$-fixed, this implies  that for every edge $e \in E(H)$, the vertex $v^*(e)$ exists. Further, if $e$ and $e'$ are distinct edges, then $v^*(e) \ne v^*(e')$. This in turn implies that for every vertex $u \in V(H)$ there is a unique edge, $e^*(u)$, such that $ v^*(e^*(u))=u$.

Let $V_1$ be the set of all vertices of color~$1$ in $\cC$ and let $V_2$ be the set of all vertices of color~$2$ in $\cC$. For each vertex $u \in V_1$, the edge $e^*(u)$ contains three vertices in $V_2$, while for each vertex $v \in V_2$, the edge $e^*(v)$ contains one vertex in $V_2$. Thus, the sum of the degrees of the vertices in $V_2$ is $3|V_1| + |V_2|$, implying that the average degree of a vertex in $V_2$ is $(3|V_1|+|V_2|)/|V_2|$. Since $H$ is $4$-regular, this value has to be~$4$, which implies that
$|V_1|=|V_2|$.

Let $H_1^*$  be the hypergraph with vertex set $V(H_1^*) = V_1$ and with edge set defined as follows: for every vertex $u \in V_2$ add the edge $e_u = (e^*(u) \setminus \{u\})$ to $H_1^*$. We note that each vertex $v \in V_1$ belongs to one edge $e^*(v)$ of $H$ and to three edges of the type $e^*(u)$ where $u \in V_2$.
Thus, by construction, $H_1^*$ is a $3$-regular $3$-uniform hypergraph.
Analogously, we define $H_2^*$ be the hypergraph with vertex set $V(H_2^*) = V_2$ and with edge set defined as follows:  for every $u \in V_1$ add the edge $e_v = (e^*(v) \setminus \{v\})$ to $H_2^*$. By construction, $H_2^*$ is a $3$-regular $3$-uniform hypergraph.
Let $C_1^1, \ldots, C_{\ell_1}^1$ be the components of $H_1^*$ where $\ell_1 \ge 1$, and let $C_1^2, \ldots, C_{\ell_2}^2$ be the components of $H_2^*$ where $\ell_2 \ge 1$. Let $i_1=1$.

Let $u_1 \in V(C_{i_1}^1)$ and let $i_2$ be defined such that $e_{u_1}$ is an edge in $C_{i_2}^2$. We note that $C_{i_2}^2 - e_{u_1}$ contains at most three components. Further, every component of $C_{i_2}^2 - e_{u_1}$ has fewer edges than vertices as the degrees of its vertices are at most~$3$ and it contains a vertex of degree at most~$2$, namely a vertex contained in the deleted edge $e_{u_1}$.
Therefore applying Lemma~\ref{lem1} to each component of $C_{i_2}^2 - e_{u_1}$, we obtain a $2$-coloring of the component that contains a free vertex. Combining these $2$-colorings in each component, produces a $2$-coloring of $C_{i_2}^2 - e_{u_1}$ that contains at least one free vertex. Let $u_2$ be a free vertex of $C_{i_2}^2 - e_{u_1}$.

Let $i_3$ be defined such that $e_{u_2}$ is an edge in $C_{i_3}^1$. As before, applying  Lemma~\ref{lem1} to each component of $C_{i_3}^1 - e_{u_2}$, we obtain a $2$-coloring with a free vertex, say $u_3$.

Continuing the above process we obtain a sequence $i_1, i_2, i_3, i_4, \ldots$. As there are finitely many components of $H_1^*$ and $H_2^*$, we note that there must exist integers $\ell$ and $k$, such that $i_{\ell}, i_{\ell+2}, i_{\ell+4}, \ldots, i_{k-2}$ are all distinct and
$i_{\ell+1}, i_{\ell+3}, i_{\ell+5}, \ldots, i_{k-1}$ are all distinct and
$i_\ell = i_k$.

Renaming components if necessary, we may assume we had started with $u_\ell$ instead of $u_1$ and we may assume $\ell=1$ and $i_1=1$, $i_2=1$, $i_3=2$, $i_4=2$, $i_5=3$, $\ldots$, $i_{k-2}=(k-1)/2$, $i_{k-1} = (k-1)/2$ and $i_k=1$. By Lemma~\ref{lem1}, the hypergraph $C_{i_k}^1 - e_{u_{k-1}} = C_1^1 - e_{u_{k-1}}$ contains at least two free vertices. We may without loss of generality assume that the free vertex $u_k$ of $C_1^1 - e_{u_{k-1}}$ was chosen to be different from $u_1$. We now define the sets $V'$ and $E'$ by $V' = \{u_2,\ldots,u_{k-1},u_k\}$ and $E' = \{e_{u_1},e_{u_2}, \ldots, e_{u_{k-1}} \}$. Further we define
\[
V^* = \bigcup_{i=1}^{\frac{k-1}{2}} V(C_i^1) \cup  V(C_i^2)
\hspace*{0.75cm} \mbox{and} \hspace*{0.75cm}
E^* = \bigcup_{i=1}^{\frac{k-1}{2}} E(C_i^1) \cup  E(C_i^2).
\]

If $e_u$ is an edge in $H_1^*$ or $H_2^*$, then let $(e_u)^H$ be the $4$-edge containing the vertices $V(e_u) \cup \{u\}$; that is, $(e_u)^H$ is the original $4$-edge in $H$ that gave rise to $e_u$. We now define
\[
E^{**} = \{ (e)^{H} \, | \, e \in E^* \} \hspace*{0.75cm} \mbox{and} \hspace*{0.75cm}
E'' = \{ (e)^{H} \, | \, e \in E' \}.
\]

We note that $E' \subset E^*$ and every edge in $E^*$ is a $3$-edge. Further, $E'' \subset E^{**}$ and every edge in $E^{**}$ is a $4$-edge. Considering the $2$-colorings we obtained above we can $2$-color the vertices of $V^* \setminus V'$ such that all edges of $E^* \setminus E'$ (and therefore also all edges of $E^{**} \setminus E''$) contain vertices of both colors. Interchanging the colors of all vertices in $C_1^1$ if necessary (by recoloring vertices of color~$i$ with color $3-i$ for $i \in [2]$ in the original $2$-coloring of $C_1^1$), we may assume that $(e_{u_1})^H$ also contains vertices of both colors, noting that $u_1 \ne u_k$. Coloring the vertex $u_i$ we can make sure that the edge $(e_{u_i})^H$ contains vertices of both colors, for each $i \in [k-1] \setminus \{1\}$.
We have now $2$-colored all vertices in $V^*$ except for the vertex $u_k$ (which is still uncolored) such that all edges in $E^{**}$ contain vertices of both colors.  Let $\cC^*$ denote the resulting $2$-coloring of the vertices of $V^*$. If $V^* = V(H)$, then $\cC^*$ is a $2$-coloring of $H$ with a free vertex, and we are therefore done and the desired result follows. Hence, we may assume that $V^* \ne V(H)$.

If the $|V^*|$ edges in $\{ e_v \, | \, v \in V^*\}$ are exactly the edges in $E^*$, then $H$ would not be connected, a contradiction. Therefore, there must be a vertex $v \in V^*$ where $e_v \not\in E^*$. Hence, $e_v \in C_a^b$ for some $a>(k-1)/2$ and $b \in [2]$. Applying  Lemma~\ref{lem1} to each component of $C_a^b - e_{v}$, we obtain analogously as before a $2$-coloring of $C_a^b - e_{v}$ with a free vertex, say $x_a^b$. If the vertex $v$ was already colored in $\cC^*$, then interchange all colors in the $2$-coloring $\cC^*$, if necessary, in order to guarantee that the edge $(e_v)^H$ contains vertices of both colors. If $v$ was not colored, then color it such that $(e_v)^H$ contains  vertices of both colors. In both cases, we produce a $2$-coloring of the vertices of $V^* \cup V(C_a^b)$ with the vertex $x_a^b$ uncolored (and possibly also the vertex $u_k$ uncolored if $v \ne u_k$) such that all edges in $E^{**} \cup \{ (e)^H \, | \, e \in E(C_a^b) \}$ contain vertices of both colors.

Repeating the above process (with the new $V^*$ being $V^* \cup V(C_a^b)$ and the new $E^{**}$ being $E^{**} \cup \{ (e)^H \, | \, e \in E(C_a^b) \}$), we will eventually have $V^* = V(H)$ and produce a $2$-coloring of $H$ with a free vertex. This completes the proof of Theorem~\ref{t:main1}.~\qed

\section{Closing Remarks}

In this paper, we establish a surprising connection between NAE-$3$-SAT and $2$-coloring of hypergraphs. We prove that every connected non-trivial instance of NAE-$3$-SAT with maximum degree~$3$ is nae-satisfiable (and contains a free variable) if the number of clauses is less than the number of variables. Using this property, we strength a beautiful result due to Carsten Thomassen~\cite{Th92} that every $4$-regular $4$-uniform hypergraph is $2$-colorable, which itself is a strengthening of a powerful result due to Alon and Bregman~\cite{AlBru88}.

As remarked earlier, our result (see Theorem~\ref{t:main1}) is best possible in the sense that there exist $4$-regular $4$-uniform hypergraphs with only free vertex; that is, every free set in such a hypergraph has size~$1$. We believe that every connected $4$-regular $4$-uniform hypergraph with sufficiently large order contains a free set of size~$2$. Due to the complement of the Fano plane the order of such a hypergraph is more than~$7$. It is possible that every connected $4$-regular $4$-uniform hypergraph of order at least~$8$ contains a free set of size~$2$. We further suspect that every connected $4$-regular $4$-uniform hypergraph with sufficiently large order~$n$ contains a free set of size $C \times n$, where $C > 0$ is some constant.

As remarked previously, we have subsequently used our NAE-$3$-SAT property, given by Theorem~\ref{main_nae_3sat}, to solve other seemingly unrelated hypergraph conjectures (such as a conjecture on the fractional disjoint transversal number) that seem difficult to solve using a purely hypergraph approach.

An interesting line of future research would be to determine for larger values of~$k \ge 4$ which connected non-trivial NAE-$k$-SAT instances are nae-satisfiable given the maximum degree, number of variables and number of clauses, and to apply such results to solve open problems and conjectures related to $k$-uniform hypergraphs. We believe this  would be a very interesting avenue of research to explore.

\newpage

\end{document}